\documentclass[letterpaper, 10 pt, conference]{ieeeconf}  
\makeatletter
\let\NAT@parse\undefined
\makeatother

\IEEEoverridecommandlockouts                              

\usepackage{graphicx}   
\usepackage{comment} 
\usepackage{amssymb,amsmath, mathrsfs,mathtools}
\usepackage[font=footnotesize,labelfont=bf]{caption}
\usepackage{setspace}
\usepackage{tabularray}
\usepackage{multirow}
\usepackage{array}
\usepackage{bigints}
\usepackage{xurl}
\usepackage{cite}
\usepackage{graphicx}   
\usepackage{longtable}
\usepackage{xurl}
\usepackage{authblk}
\usepackage{lineno}
\usepackage{threeparttable}
\usepackage[colorlinks,allcolors=blue]{hyperref}
\newtheorem{theorem}{Theorem}
\newtheorem{definition}{Definition}

\title{\LARGE \bf Highly Efficient Optimal Control for Lyophilization via Simulation of Discrete/Continuous Mixed-index Differential-algebraic Equations}

\author{Prakitr Srisuma and Richard D. Braatz 
\thanks{The authors are with the Massachusetts Institute of Technology, Cambridge, MA 02139. Email:
        \{prakitrs, braatz\}@mit.edu}%
}

\begin{document}

\maketitle
\thispagestyle{empty}
\pagestyle{empty}

\begin{abstract}

This article presents a highly efficient optimal control algorithm and policies for lyophilization (also known as freeze drying). The optimal solutions and control policies are derived using an extended version of the simulation-based algorithm, which reformulates the optimal control problem as a hybrid discrete/continuous system of mixed-index differential-algebraic equations and subsequently calculates the optimal control vector via simulation of the resulting DAEs. Our algorithm and control policies are demonstrated via a number of case studies that encompass various lyophilization and optimal control strategies. All the case studies can be solved within less than a second on a normal laptop, regardless of their complexity. The method is several orders of magnitude faster than the traditional optimization-based techniques while giving similar/better accuracy. The proposed algorithm offers an efficient and reliable framework for optimal control of lyophilization, which can also be extended to other similar systems with phase transitions.

\end{abstract}

\section{Introduction} \label{sec:Intro}

Lyophilization (aka freeze drying) is a crucial process in the (bio)pharmaceutical industry used to improve the stability of drug products \cite{Fissore2018Review}, including its recent application to mRNA vaccines for COVID-19. Advances in lyophilization technology play an important role in enhancing storage and distribution of drug products and the entire biopharmaceutical manufacturing in general.

Various control strategies have been studied for lyophilization; we refer to \cite{Daraoui2010ReviewControl,Fissore2019Book} for a brief summary of those studies. Optimal control of lyophilization has also been explored to a certain extent, with some promising results and important observations. For example, \cite{Litchfield1982PressureControl} found that manipulating the heat input resulted in the largest reduction in drying times, while varying the pressure produced only a small change. In \cite{Sadikoglu1998Indirect}, the optimal control policies were obtained using the indirect method (i.e., via variational calculus) based on a one-dimensional physics-based model, resulting in a 40\% decrease in drying time. The method was later extended to a multidimensional model \cite{Sadikoglu2003Indirect}. The direct method, which relies on control vector parameterization, has also been applied for optimal control of lyophilization \cite{Pisano2010Control,Pisano2011MPC,Antelo2012Direct,Vilas2020Direct}. In \cite{Vilas2020Direct}, a 30\% drying time reduction was reported.

One of the typical challenges in optimal control problems is its high computational cost associated with solving large-scale nonlinear optimization problems \cite{Pisano2010Control,Pisano2011MPC,Bano2020LumpedDrying}. Solving an optimal control problem online (MPC) or with probabilistic uncertainty (stochastic control) is even more computationally expensive, limiting their applications in real-world manufacturing systems. Recently, a new class of optimal control algorithms, the simulation-based method, was proposed \cite{Berliner2022battery,Srisuma2025SimDAE}. This method transforms an optimal control problem into a system of differential-algebraic equations (DAEs), in which the optimal control and state trajectories can be obtained by simulation of the resulting DAEs, without solving any optimization. This method was shown to be several orders of magnitude faster than the optimization-based approaches while maintaining similar/better accuracy \cite{Srisuma2025SimDAE}. This algorithm, however, is currently limited to optimal problems reformulated as index-1 DAEs \cite{Berliner2022battery} or high-index DAEs with smooth trajectories in the absence of path constraints \cite{Srisuma2025SimDAE}.

This article presents a highly efficient optimal control algorithm and policies for lyophilization. The main contributions of this work are twofold. On the algorithm side, we develop an extended version of the simulation-based method for handling a hybrid discrete/continuous system of mixed-index DAEs, which expands the application of the simulation-based method to a broader class of optimal control problems, including those involving index-1 DAEs, high-index DAEs, non-smooth control trajectories, and path constraints. On the application side, we discuss three control policies for lyophilization and derive an equivalent system of DAEs corresponding to each policy. The proposed control policies are then implemented via the simulation-based method to solve various optimal control case studies associated with lyophilization.

This article is organized as follows. Section \ref{sec:Model} describes the lyophilization system and summarizes the important model equations. Section \ref{sec:OCP_algo} discusses the simulation-based method, including its theory, extension, algorithm, and implementation. Section \ref{sec:OCP_policy} summarizes the three control policies for lyophilization and derives a system of DAEs corresponding to each policy. Section \ref{sec:Results} demonstrates the extended simulation-based method via several case studies on optimal control of lyophilization. Finally, Section \ref{sec:Conclusion} summarizes the study and suggests some future directions.


\section{Mechanistic Modeling} \label{sec:Model}
\subsection{Process description}
This study considers continuous lyophilization of suspended vials, the cutting-edge continuous lyophilization technology recently proposed by \cite{Capozzi2019ContLyo_SuspendedVials}. With this technology, a number of vials are suspended and move continuously through the process without any contact between the vials and shelf/chamber. For optimal control, we consider the primary drying step only as it is the most time-consuming and expensive step in lyophilization, and hence the main target for improvement and optimization. The model used in this work is mainly based on the state-of-the-art lyophilization model proposed by \cite{Srisuma2025ContModel}, in which the key model equations are summarized below. We refer to \cite{Srisuma2025ContModel} for the detailed derivation of all equations and model validation.
\vspace{-5pt}
\begin{figure}[ht]
    \centering
    \includegraphics[scale=1.1]{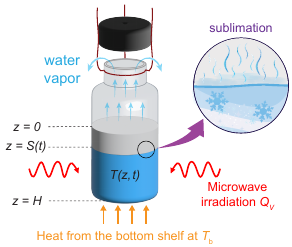}
    \caption{Schematic diagram showing the primary drying model.} 
    \label{fig:Schematic}
\end{figure}
\vspace{-5pt}
The model for primary drying is formulated in the rectangular coordinate system with one spatial dimension ($z$) and time ($t$) (Fig.\ref{fig:Schematic}). The product is heated by the below shelf, with the shelf temperature $T_\textrm{b}$. By assuming that the supplied heat is used in the frozen region only, the energy balance for the frozen region can be described by the partial differential equation (PDE)
\begin{equation} \label{eq:1st_energy}
    \rho_\textrm{f} C_{p,\textrm{f}}\dfrac{\partial T}{\partial t} = k_\textrm{f}\dfrac{\partial^2 T}{\partial z^2}+ \frac{Q_\textrm{rad}}{V_\textrm{f}}, \qquad  S < z < H, 
\end{equation}
where $T(z, t)$ is the temperature, $S(t)$ is the sublimation front/interface position, $k$ is the thermal conductivity, $\rho$ is the density, $C_p$ is the heat capacity, $V$ is the volume, $H$ is the height of the product, and the subscript `f' denotes the frozen region. The radiative heat transfer from the sidewall $Q_\textrm{rad}$ is
\begin{equation} \label{eq:1st_rad}
    Q_\textrm{rad} = \sigma A_r \mathcal{F}_\textrm{side} \!\left(T_\textrm{c}^4-T^4\right),
\end{equation}
where $T_\textrm{c}$ is the chamber wall temperature, $A_r=\pi dH$ is the side area of the product, $\mathcal{F}_\textrm{side}$ is the transfer factor of the side surface, and $\sigma$ is the Stefan-Boltzmann constant. The mass balance of water at the sublimation front can be described by the ordinary differential equation (ODE)
\begin{equation}\label{eq:1st_interface}
    \frac{dS}{dt} = \frac{N_\textrm{w}}{\rho_\textrm{f}-\rho_\textrm{e}},
\end{equation}
where $N_\textrm{w}$ is the sublimation flux and $\rho_\textrm{e}$ is the effective density of the dried region above the sublimation front. The driving force for mass transfer at the sublimation interface is \cite{Pikal2005Model,Fissore2018Review,Bano2020LumpedDrying}
\begin{equation}\label{eq:1st_flux}
    N_\textrm{w} = \frac{p_\textrm{w,sat}-p_\textrm{w,c}}{R_\textrm{p}},
\end{equation}
where $p_\textrm{w,sat}$ is the saturation/equilibrium pressure of water, $p_\textrm{w,c}$ is the partial pressure of water in the chamber (environment), and $R_\textrm{p}$ is the mass transfer resistance.

With the proper boundary conditions and numerical methods (see \cite{Srisuma2025ContModel}), the above model equations can be written as a system of ODEs 
\begingroup
\allowdisplaybreaks
\begin{align}
&\dfrac{dT_1}{dt} = f_1(T_1, T_2, S),  \label{eq:model1}\\
&\dfrac{dT_i}{dt} = f_2(T_{i-1}, T_i, T_{i+1}, S), \quad  i = 2, \dots, n-1,  \label{eq:model2}\\
&\dfrac{dT_n}{dt} = f_3(T_{n-1}, T_n, S, T_\textrm{b}),  \label{eq:model3}\\
&\dfrac{dS}{dt} = f_4(T_1, S), \label{eq:model4}
\end{align}
\endgroup
where $f_1, f_2, f_3, f_4$ are the nonlinear functions, $n$ is the number of grid points for spatial discretization, $T_i$ is the discretized product temperature; $i=1$ denotes the top interface (aka sublimation interface); and $i=n$ denotes the bottom of the product. The initial conditions are
\begin{align}
    T(z, t_0) &= T_0, \qquad 0\leq z \leq H, \label{eq:1st_iniT}  \\
    S(t_0) &=  0, \label{eq:1st_iniS} 
\end{align}
where $t_0$ is the initial time. The primary drying model is simulated until the interface position is equal to the height of the product, i.e., $S = H$, indicating that there is no frozen material left, which marks the end of the primary drying step. The final time is denoted as $t_f$.

\section{Simulation-based Method} \label{sec:OCP_algo}
The optimal control formulation for continuous lyophilization is
\begingroup
\allowdisplaybreaks
\begin{align}\label{eq:OptControlGeneral}
&\begin{aligned}
\min_{T_\textrm{b}(t)}  \   &\mathcal{M}(T(t_f),S(t_f)) + \int_{t_0}^{t_f}\mathcal{L}(T(t)S(t),T_\textrm{b}(t),t) dt 
\end{aligned} \\
& \: \textrm{s.t.  Equations \eqref{eq:1st_iniT}, \eqref{eq:1st_iniS}, and \eqref{eq:model1}--\eqref{eq:model4}} , \nonumber \\
&\quad \ \ \ \, h(T(t),S(t),t) \leq 0, \label{eq:path} \\
&\quad \ \ \  \, T_{\textrm{b,min}} \leq T_\textrm{b}(t) \leq T_{\textrm{b,max}} \label{eq:bound},
\end{align}
\endgroup
where \eqref{eq:path} is the path constraints and \eqref{eq:bound} are the bounds on a control vector. 

Traditional numerical algorithms for solving optimal control problems, including direct and indirect methods, have been discussed extensively in the literature \cite{Rodrigues2014OCoverview,Nolasco2020OptimalControl,Srisuma2025SimDAE}. To summarize, direct methods, which entail control vector parameterization and nonlinear programming, are mostly used nowadays due to their ease of implementation and applicability to a wide range of problems. Nevertheless, the key limitation of direct methods is their high computational cost, especially when solving large-scale nonlinear optimization. For example, a computation time of 4.5 h was reported in \cite{Bano2020LumpedDrying}, which is clearly not feasible for real-time applications in manufacturing. In addition, the performance and accuracy of direct methods can vary significantly depending on several factors, e.g., initial guesses, control vector parameterization, selected optimization solvers, and optimality tolerances. Note that, since these techniques require an optimization solver, we refer to it as the \textit{optimization-based} methods. 

\subsection{Overview}
A new class of optimal control algorithms, the \textit{simulation-based} method, has been recently proposed \cite{Berliner2022battery,Srisuma2025SimDAE}. The simulation-based method reformulates an optimal control problem as a system of DAEs, in which the optimal control trajectory can be obtained via \textit{simulation} of the resulting DAEs using a proper DAE solver instead of solving an \textit{optimization} problem. Since the problem is solved via simulation rather than optimization, the computation is much faster by several orders of magnitude \cite{Srisuma2025SimDAE}. 

The simulation-based method was first demonstrated for the fast charging of lithium-ion batteries \cite{Berliner2022battery}, in which the optimal control problem was reformulated as a mixed continuous-discrete system of index-1 DAEs, which was implemented in Julia.

In many engineering applications, the reformulation can result in high-index DAEs, which are not compatible with solvers designed for index-1 systems. In \cite{Srisuma2025SimDAE}, the simulation-based method was formalized and generalized to high-index DAEs (up to index-20), in which the method was shown to be several orders of magnitude faster than all the optimization-based methods while giving similar/better accuracy, regardless of the differential index. The capability to handle high-index DAEs significantly broadens the applications of the simulation-based method to a wider class of optimal control problems. Nevertheless, \cite{Srisuma2025SimDAE} primarily focuses on case studies with smooth control trajectories in the absence of path constraints.

The algorithm proposed in this article is a major extension of both \cite{Berliner2022battery} and \cite{Srisuma2025SimDAE}, where an optimal control problem is reformulated as a hybrid discrete/continuous system of mixed-index DAEs, enabling the method to solve problems with path constraints and/or non-smooth control trajectories. This extension greatly expands the application of the simulation-based method to more optimal control problems.

\subsection{Theory and algorithm}

\begin{definition}[Mixed-index DAEs]
In this work, mixed-index DAEs refer to a system of DAEs that consists of both index-1 and high-index DAEs, in which the differential index can vary over time. Index-1 DAEs refer to DAEs with a differential index of 1. High-index DAEs refer to DAEs with a differential index greater than 1.

\end{definition}

\begin{theorem} \label{theorem1}
Some optimal control problems in the form of \eqref{eq:OptControlGeneral} can be reformulated as a hybrid discrete/continuous system of mixed-index DAEs
\begin{equation} \label{eq:DAE}
\begin{gathered}
    g_1(T(t),S(t),T_\textrm{b}(t),t) = 0,\\
    g_2(T(t),S(t),T_\textrm{b}(t),t) = 0,\\
    \vdots  \\ 
    g_m(T(t),S(t),T_\textrm{b}(t),t) = 0,
\end{gathered}    
\end{equation} 

where $m$ is the total number of DAE systems resulting from the reformulation, which varies among problems. Consequently, the optimal solution to the original problem \eqref{eq:OptControlGeneral} can be obtained by solving \eqref{eq:DAE} instead.
\end{theorem}

The simulation-based algorithm consists of two main steps. The first step is to reformulate the optimal control problem as an equivalent system of DAEs as shown by \eqref{eq:DAE}. This reformulation could rely on the mechanistic understanding of the system, optimality conditions, or analysis of the problem structure; we show some examples in the case studies presented in Sections \ref{sec:OCP_policy} and \ref{sec:Results}. The second step is to solve the resulting DAEs \eqref{eq:DAE} correctly.

The key idea of the simulation-based method is to change fron an optimization problem to a simulation (DAE) problem. Consequently, the control trajectory can be obtained by numerically solving (i.e., simulating) the resulting DAEs using a proper DAE solver. This approach has three main benefits. First, the algorithm is highly efficient due to the absence of numerical optimization. Second, the accuracy of the solution and computational performance are more consistent than the optimization-based methods because the method does not rely on initial guesses, optimization tolerances, or control vector parameterization. Finally, the implementation of the algorithm requires mechanistic understanding about the process, and so the results are highly interpretable.

\subsection{Implementation}  \label{sec:Implementation}
Since this work addresses the most complex and general version of the simulation-based method to date, the algorithm and implementation is significantly more complicated than the previous iterations of this algorithm \cite{Berliner2022battery,Srisuma2025SimDAE}. 

For optimal control problems reformulated as an index-1 system, the resulting DAEs can be solved easily with any DAE solvers, without any further technique. In this work, MATLAB's \texttt{ode15s} is used. MATLAB's \texttt{ode15s} is highly efficient due to its adaptive time-stepping scheme. For high-index DAE system, the best option, as described in \cite{Srisuma2025SimDAE}, is to use the DAE solver developed as part of GEEKO \cite{GEKKO}, a Python package for optimization and machine learning. Unlike MATLAB's \texttt{ode15s}, GEKKO’s solvers do not support adaptive time-stepping, which results in significantly slower computation. As of now, this limitation represents the primary bottleneck in the computational efficiency of the simulation-based method. Nevertheless, GEKKO remains the most suitable tool currently available for solving such high-index systems.

In \eqref{eq:DAE}, the simulation-based method entails selecting and solving a system of DAEs corresponding to the optimal control policy over time ($g_1, g_2, ..., g_m$). Each DAE system $g_i$ represents one unique control policy. When the optimal control policy changes from one to another, e.g., from $g_1$ to $g_2$, we denote it as an \textit{event}. During simulation, continuous event detection is required to ensure that the optimal control policy is always selected and simulated. When an event is detected, meaning that the current policy is not optimal, a system of DAEs is reinitialized to represent the new optimal policy (i.e., policy switching). Repeat these steps until the termination criterion is met. 

Continuous event detection during the simulation of DAEs is dependent on the choice of a DAE solver. For example, MATLAB and Julia have their own built-in functions for event handling; these functions are known as event functions in MATLAB and callback functions in Julia. However, this study deals with mixed-index DAEs, which requires both MATLAB's \texttt{ode15s} and GEKKO's DAE solver, and so event detection and policy switching require efficient and seamless communication between both solvers. GEKKO's DAE solver, however, does not support any event functions, and so event detection during the simulation of high-index DAEs is more complicated. To address this limitation, we propose to use GEKKO for solving high-index DAEs and then return the solution to MATLAB for event detection and policy switching as follows. First, simulate a system of high-index DAEs using GEKKO under the assumption that no events occur within the time horizon of interest. Then, transfer the control trajectory obtained from GEKKO to MATLAB. Finally, simulate the model equations \eqref{eq:model1}--\eqref{eq:1st_iniS} using MATLAB's \texttt{ode15s} for event detection, with the control trajectory from GEKKO specified as an input \footnote{As the control/input is now specified, the resulting DAE system has a differential index of 1, which can be now solved using \texttt{ode15s}.}.

All simulations were performed in MATLAB R2024b, with GEKKO called and executed in Python 3.10, on a computer equipped with an AMD Ryzen™ 9 5900HS processor (8 cores) and 32 GB RAM running 64-bit Windows 11.

\section{Control Policies for Lyophilization} \label{sec:OCP_policy}
This section describes three different control policies for lyophilization, the reformulation technique, and the final DAE system corresponding to each control policy.

\subsection{Policy 1: Maximum heat input}
The first policy is the simplest and most typical that is widely used lyophilization processes. This policy focuses on maximizing the total heat input to accelerate the drying process. This policy corresponds to setting $T_\textrm{b}(t)$ to its upper limit, resulting in the DAE system $g_1$:
\begingroup
\allowdisplaybreaks
\begin{equation} \label{eq:policy1}
\begin{aligned}
     & T_\textrm{b}(t) = T_\textrm{b,max}, \\
    &\textrm{Equations \eqref{eq:1st_iniT}, \eqref{eq:1st_iniS}, and \eqref{eq:model1}--\eqref{eq:model4}}, 
\end{aligned}
\end{equation}
\endgroup
which is an index-1 DAE system. Hence, the DAEs correspond to Policy 1 can be solved using MATLAB's \texttt{ode15s}. 
In lyophilization, Policy 1 is set as the optimal policy by default provided that there is no other active constraint. With no other active constraint, it implies that the system is operated in the design space, and so the maximum heat input would lead to fastest drying.

\subsection{Policy 2: Product temperature tracking}
Policy 2 focuses on manipulating the shelf temperature such that the product temperature is kept at a desired setpoint $T_\textrm{sp}$. In the context of lyophilization, there are various temperature limits that should be considered, e.g., collapse temperature, glass transition temperature, and melting temperature \cite{Fissore2018Review}. It is important to ensure that the product temperature does not exceed those limits. This policy corresponds to the optimal control problem
\begingroup
\allowdisplaybreaks
\begin{equation}\label{eq:policy2_original}
\begin{aligned}
&\min_{T_\textrm{b}(t)} \ \int_{t_0}^{t_f}\left(T_n(t)-T_\textrm{sp}\right)^2 dt \\
& \: \textrm{s.t. \ Equations \eqref{eq:1st_iniT}, \eqref{eq:1st_iniS}, and \eqref{eq:model1}--\eqref{eq:model4}},
\end{aligned} 
\end{equation}
\endgroup
where $T_n$ (product temperature at the bottom, $i=n$) is used to represent the product temperature because, according to the typical design of lyophilizers, the product temperature tends to be highest at the bottom due to temperature gradients resulting from heat transfer. 

To implement the simulation-based method, \eqref{eq:policy2_original} needs to be reformulated as a system of DAEs. This objective function is minimized when the rate of temperature change is equal to the setpoint. Replacing the objective function in \eqref{eq:policy2_original} with the algebraic equation $T(t) = T_{\textrm{sp}}$ results in the DAE system $g_2$:
\begingroup
\allowdisplaybreaks
\begin{equation} \label{eq:policy2}
\begin{aligned}
     & T_n(t) = T_\textrm{sp}, \\
     & T_\textrm{b}(t_0) = T_\textrm{b0}, \\
    &\textrm{Equations \eqref{eq:1st_iniT}, \eqref{eq:1st_iniS}, and \eqref{eq:model1}--\eqref{eq:model4}},
\end{aligned}
\end{equation}
\endgroup
which is an index-2 DAE system, hence GEKKO's DAE solver.

\subsection{Policy 3: Sublimation flux tracking}
Policy 3 focuses on manipulating the shelf temperature such that the sublimation flux is kept at a desired setpoint to avoid vapor accumulation in the drying chamber \cite{Fissore2018Review} and control the interface position and drying time \cite{Daraoui2010ReviewControl}. 

Since the sublimation flux $N_\textrm{w}$ is directly proportional to the interface velocity $dS/dt$ (see \eqref{eq:1st_interface}), controlling either process variable yields the same results. In this study, we focus on controlling the interface velocity as it can be easily related to the interface position $S$, one of the critical process variables in primary drying. As such, Policy 3 can be represented by the optimal control problem 
\begingroup
\allowdisplaybreaks
\begin{equation}\label{eq:policy3_original}
\begin{aligned}
&\min_{T_\textrm{b}(t)} \  \int_{t_0}^{t_f}\left(\frac{dS}{dt}-v_\textrm{sp}\right)^2 dt \\
& \: \textrm{s.t. \ Equations \eqref{eq:1st_iniT}, \eqref{eq:1st_iniS}, and \eqref{eq:model1}--\eqref{eq:model4}},
\end{aligned}
\end{equation}
\endgroup
where $v_\text{sp}$ is the target interface velocity. Similar to Policy 2, the bound on $T_\textrm{b}(t)$ is not included in \eqref{eq:policy3_original} as it is already considered in Policy 1.

It is obvious that the objective function of \eqref{eq:policy3_original} is minimized when $dS/dt$ is equal to $v_\text{sp}$. Hence, we enforce this condition by modifying the model equation that describes the interface velocity \eqref{eq:model4}, resulting in the DAE system $g_3$:
\begingroup
\allowdisplaybreaks
\begin{equation} \label{eq:policy3}
\begin{aligned}
     & \frac{dS}{dt} = v_\textrm{sp}, \\
     & f_4(T,S) = v_\textrm{sp}, \\
     & T_\textrm{b}(t_0) = T_\textrm{b0}, \\
    &\textrm{Equations \eqref{eq:model1}--\eqref{eq:model3}, \eqref{eq:1st_iniT}, and \eqref{eq:1st_iniS}}.
\end{aligned}
\end{equation}
\endgroup
The differential index of \eqref{eq:policy3} is $n+1$, which is dependent on spatial discretization of the PDE. If the spatial discretization of the domain is made finer to increase numerical accuracy, the differential index of \eqref{eq:policy3} increases. In this work, $n$ is set to 20, and hence \eqref{eq:policy3} is an index-21 DAE system, which can be solved using GEKKO DAE's solver. 

\subsection{Overall policy selection}
To summarize, there are three control policies considered in this study, hence three DAE systems: $g_1,g_2,g_3$. By default, for situations where no event is detected (i.e., no active inequality constraint), Policy 1 ($g_1$) is always selected. If an event is detected, the simulation-based method will switch to that policy. For example, if the sublimation flux reaches its upper limit, Policy 3 ($g_3$) will be selected. Every control policy is represented by a unique system of DAEs, which needs to be solved using a proper DAE solver.

\begin{figure*}[ht]
    \centering
    \includegraphics[scale=.9]{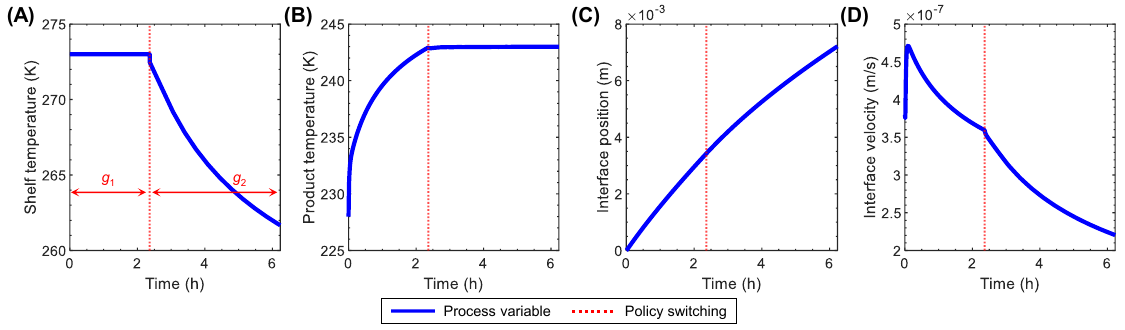}
    \caption{Optimal solution to Problem 1 showing the time evolution of the (A) shelf temperature, (B) maximum product temperature, (C) interface position, and (D) interface velocity. The optimal control policy switches from Policy 1 to Policy 2 at 2.4 h.} 
    \label{fig:OCP_Case2}
\end{figure*}


\section{Results and Discussion} \label{sec:Results}
This section implements the extended simulation-based method to solve various optimal control problems for lyophilization. All problems were solved using the simulation-based method, in which the initialization, event detection, policy selection/switching, and termination were executed automatically as described in Section \ref{sec:Implementation}. The default model parameters are available in the software released with this work (see Data Availability), while problem-specific parameters (e.g., setpoints, bounds) are reported in the corresponding case study.

\subsection{Problem 1: Minimizing the drying time}
Problem 1 represents the most common optimal control scenario that can be found and implemented practically, e.g., as demonstrated in \cite{Sadikoglu1998Indirect,Pisano2010Control,Antelo2012Direct}. This case study focuses on finding the optimal trajectory of the shelf temperature that minimizes the drying time while satisfying the stability constraints associated with product temperature, corresponding to the optimal control problem
\begin{equation}
\begin{aligned}\label{eq:OptControlCase2}
&\min_{T_\textrm{b}(t)} \   t_f \\
& \: \textrm{s.t. \ Equations \eqref{eq:1st_iniT}, \eqref{eq:1st_iniS}, and \eqref{eq:model1}--\eqref{eq:model4}}, \\
&\quad \ \ \ \ \, T(z,t) \leq 243 \ \textrm{K},   \\ 
&\quad \ \ \ \ \, 228 \ \textrm{K} \leq T_\textrm{b}(t) \leq 273  \ \textrm{K}, 
\end{aligned} 
\end{equation}
where the upper limit on the product temperature is set to 243 K, and the shelf temperature can be varied between 228 K and 273 K.

With the simulation-based method, \eqref{eq:OptControlCase2} can be solved in about $0.58\pm0.01$ s. The optimal solution consists of Policies 1 and 2 (Fig.~\ref{fig:OCP_Case2}), which is consistent with the problem formulation \eqref{eq:OptControlCase2}. The shelf temperature is held constant at its upper bound of 273 K (Policy 1) since the initial time until 2.4 h (Fig.~\ref{fig:OCP_Case2}A). At 2.4 h, the maximum product temperature reaches its upper limit of 243 K (Fig.~\ref{fig:OCP_Case2}B), and thus the optimal control policy switches from Policy 1 to Policy 2. After 2.4 h, the shelf temperature decreases such that the maximum product temperature is maintained at the upper limit (Figs.~\ref{fig:OCP_Case2}AB). The optimal shelf temperature in this case study exhibits a similar shape and trajectory to those reported in the literature \cite{Sadikoglu1998Indirect,Pisano2010Control}. The solution from the simulation-based method also agrees with those obtained from the traditional optimization-based techniques, with the simulation based approach being several orders of magnitude faster (see the provided software and code in Data Availability for the benchmarking results).

\subsection{Problem 2: Fastest drying with multiple policies}
This problem considers a more complicated version of Problem 1, in which different constraints and setpoints are involved. Specifically, a constraint on the interface velocity is added, resulting in the optimal control problem
\begin{equation}\label{eq:OptControlCase4}
\begin{aligned}
&\min_{T_\textrm{b}(t)} \   t_f \\
& \: \textrm{s.t. \ Equations \eqref{eq:1st_iniT}, \eqref{eq:1st_iniS}, and \eqref{eq:model1}--\eqref{eq:model4}}, \\
&\quad \ \ \ \ \, T(z,t) \leq 240 \ \textrm{K},  \\
&\quad \ \ \ \ \, \frac{dS}{dt} \leq 2.8\times10^{-7} \ \textrm{m/s},  \\
&\quad \ \ \ \ \, 228 \ \textrm{K} \leq T_\textrm{b}(t) \leq 260  \ \textrm{K}.
\end{aligned} 
\end{equation}

\begin{figure*}[ht]
    \centering
    \includegraphics[scale=.9]{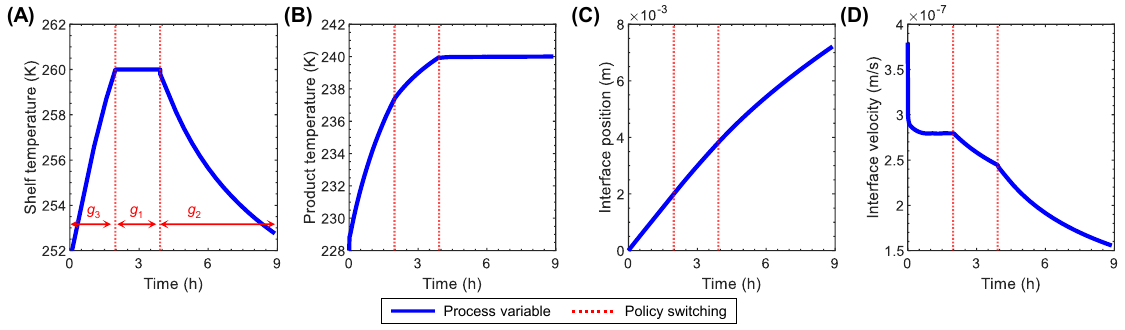}
    \caption{Optimal solution to Problem 2 showing the time evolution of the (A) shelf temperature, (B) maximum product temperature, (C) interface position, and (D) interface velocity. The optimal control policy switches from 3 to 1 at 2.0 h and from 1 to 2 at 3.9 h, respectively.} 
    \label{fig:OCP_Case4}
\end{figure*}

With the simulation-based method, \eqref{eq:OptControlCase4} can be solved in about $0.98\pm0.02$ s. In this case, the optimal solution consists of Policies 1, 2, and 3 (Fig.~\ref{fig:OCP_Case4}). At the beginning, the interface velocity exceeds its upper limit of $2.8\times10^{-7}$ m/s, and so the simulation-based method selects Policy 3 to manipulate the shelf temperature such that the interface velocity converges to $2.8\times10^{-7}$ m/s quickly (Fig.~\ref{fig:OCP_Case4}D). At $t=2$ h, the shelf temperature reaches its upper bound, and thus our algorithm switches to Policy 1 (Fig.~\ref{fig:OCP_Case4}A). The shelf temperature is maintained at its upper bound until $t=3.9$ h, when the product temperature approaches its upper limit of 240 K. Consequently, the simulation-based method switches to Policy 2 to maintain the product temperature at 240 K until the drying process is complete at about 8.9 h.

\section{Conclusion} \label{sec:Conclusion}
This article describes a highly efficient optimal control algorithm and policies for lyophilization. The proposed algorithm, the extended simulation-based method, reformulates an optimal control problem as a hybrid discrete/continuous system of mixed-index differential-algebraic equations, in which the optimal control and state trajectories can be obtained via simulation of the resulting DAEs. The approach is demonstrated for two optimal control case studies related to lyophilization. With the extended simulation-based method, all the case studies can be solved within less than a second on a normal laptop, regardless of their complexity. The method is several orders of magnitude faster than the traditional optimization-based techniques while giving similar/better accuracy. The proposed framework offers an efficient and reliable framework for optimal control of lyophilization, which can also be extended to other similar systems with phase transition.

\section*{Data Availability}  \label{sec:Code}
Software and data used in this work are available at \url{https://github.com/PrakitrSrisuma/simDAE-optimalcontrol-lyo}.

\section*{Acknowledgments} 
This research was supported by the U.S. Food and Drug Administration under the FDA BAA-22-00123 program, Award Number 75F40122C00200.

\bibliographystyle{IEEEtran}
\bibliography{IEEEabrv,reference}

\end{document}